\documentclass[12pt,a4paper]{article}
\usepackage[cp1251]{inputenc} % cp866 - DOS, cp1251 - Windows
\usepackage[russian]{babel}
\usepackage{amsfonts}

\begin{document}

\begin{center}
{\large\bf Global Solution of the Inverse Spectral Problem for
Differential Operators on a Finite Interval with Complex Weights}\\[0.2cm]
{\bf V.A.\,Yurko}
\end{center}

\thispagestyle{empty}

\noindent {\bf Abstract.} Non-self-adjoint second-order ordinary
differential operators on a finite interval with complex weights are
studied. Properties of spectral characteristics are established and
the inverse problem of recovering operators from their spectral
characteristics are investigated. For this class of nonlinear inverse
problems an algorithm for constructing the global solution is obtained.
To study this class of inverse problems, we develop ideas of the method
of spectral mappings.

\smallskip
\noindent {\bf Keywords:} differential operators, complex weight,
spectral characteristics, inverse problem, method of spectral mappings

\smallskip
\noindent {\bf Mathematics Subject Classification:} 34A55, 34B24, 47E05.

\bigskip
{\bf 1. Introduction}

\medskip
We consider the boundary value problem $L$ for the differential equation
$$
-y''(x)+q(x)y(x)=\lambda r(x)y(x),\quad 0<x<T,                   \eqno(1.1)
$$
subject to the Robin boundary conditions
$$
U(y):=y'(0)-hy(0)=0,\quad V(y):=y'(T)+Hy(T)=0,                   \eqno(1.2)
$$
and the jump conditions at an interior point $b\in(0,T)$:
$$
y(b+0)=d_1 y(b-0),\quad y'(b+0)=y(b-0)/d_1+d_2y(b-0).            \eqno(1.3)
$$
Here $\lambda$ is the spectral parameter, $q(x)$ and $r(x)$ are
complex-valued functions, $q(x)\in L(0,T),$ and $r(x)=a_k^2$ for
$x\in(b_{k-1},b_k),$ where $0=b_0<b_1=b<b_2=T.$ The numbers
$h, H, a_k$ and $d_k$ are complex, and $a_k\ne 0, d_1\ne 0.$
For definiteness, let $\arg d_1\in[0,\pi).$

We study the inverse spectral problem for the boundary value problem
(1.1)–(1.3). Inverse spectral problems consist in recovering operators
from their spectral characteristics. Such problems play an important
role in mathematics and have many applications in natural science
and technology. Inverse spectral problems are also used for solving
nonlinear integrable evolution equations of mathematical physics.
Inverse problems for the classical Sturm–Liouville operators
(when $r(x)\equiv 1,$ $d_1=1,$ and $d_2=0$) have been studied fairly
completely (see [1] and the historical review therein). Inverse problems
for arbitrary order differential operators and systems with arbitrary
characteristic numbers are more difficult. They have been solved later
by the method of spectral mappings (see the monographs [2]–[3] and the
references therein). Inverse problems on spatial networks are an important
and popular part of the inverse problem theory; in the review paper [4]
one can find the main results on inverse problems on spatial networks.
Boundary value problems with discontinuous weights and jump conditions
at interior points have been considered in many papers, but mostly
for the case with real weights. In the case when $r(x)\equiv 1$
(i.e. $a_k=1$), the boundary value problem $L$ satisfying conditions
(3) was studied in [5]–[9] and other papers. Inverse problems for a real
weight $r(x)$ were studied in [10]–[14] and other works. Inverse
problems for the boundary value problem $L$ with complex-valued weights
were studied in [15]-[16] where only uniqueness results were obtained.
Note that complex-valued weights appear, in particular, in the study of
the interaction of electromagnetic waves with layered media possessing
both dielectric and magnetic properties [17]. Moreover, a number of
problems for Sturm–Liouville equations on curves in the complex plane
can be reduced to the boundary-value problem $L$ of the form (1)–(3)
on a real interval. In the present paper, we establish properties of
the spectral characteristics for $L,$ and study the inverse spectral
problem of recovering parameters of $L$ from the given spectral
characteristics. For this class of nonlinear inverse problems
an algorithm for constructing the global solution is obtained.
To study this class of inverse problems, we develop ideas related
to the method of spectral mappings [2].

\medskip
{\bf 2. Spectral data}

\medskip
Let $l_k:=b_k-b_{k-1}$ and $a_k=r_k\exp(i\varphi_k),$ $r_k>0,$
$0\le\varphi_2<\varphi_1<\pi.$ We assume that the following regularity
condition holds: $\omega_{\pm}:=d_1a_2\pm a_1/d_1\ne 0.$ Denote by
$\Phi(x,\lambda)$ the solution of (1.1) such that (1.3) holds and
$U(\Phi)=1,$ $V(\Phi)=0.$ Let $M(\lambda):=\Phi(0,\lambda).$ We will
also use the solutions $\varphi(x,\lambda), \psi(x,\lambda),
S(x,\lambda)$ of Eq. (1.1) satisfying (1.3) and the conditions
$\varphi(0,\lambda)=1, \varphi'(0,\lambda)=h, S(0,\lambda)=0,
S'(0,\lambda)=1, \psi(T,\lambda)=1, \psi'(T,\lambda)=-H.$
Denote $D(x,\lambda,\mu):=(\lambda-\mu)^{-1}
\langle\varphi(x,\lambda),\varphi(x,\mu)\rangle,$
where $\langle y(x),z(x)\rangle:=y(x)z'(x)-y'(x)z(x).$
The function
$$
\Delta(\lambda):=\langle\varphi(x,\lambda),\psi(x,\lambda)\rangle
=-V(\varphi)=U(\psi)                                            \eqno(2.1)
$$
does not depend on $x,$ and it is called the characteristic function
for $L.$ The eigenvalues $\Lambda:=\{\lambda_k\}_{k\ge 0}$ of $L$ coincide
with the zeros of the entire function $\Delta(\lambda).$ Clearly,
$$
\Phi(x,\lambda)=S(x,\lambda)+M(\lambda)\varphi(x,\lambda)=
\psi(x,\lambda)/\Delta(\lambda),\;
M(\lambda)=\Delta_0(\lambda)/\Delta(\lambda),                  \eqno(2.2)
$$
where $\Delta_0(\lambda):=\psi(0,\lambda)=V(S).$
Using (2.1) and (2.2) one gets
$$
\langle\varphi(x,\lambda),\Phi(x,\lambda)\rangle\equiv 1.      \eqno(2.3)
$$

Let $\lambda=\rho^2, \lambda_k=\rho_k^2.$ Consider the half-planes
$\Pi_k^{\pm}:=\{\rho:\;\pm Im(\rho a_k)>0\},$ $k=1,2,$ and denote
$$
S_1=\Pi_1^{+}\cup\Pi_2^{+},\; S_2=\Pi_1^{-}\cup\Pi_2^{+},\;
S_3=\Pi_1^{-}\cup\Pi_2^{-},\; S_4=\Pi_1^{+}\cup\Pi_2^{-}.
$$
Then $S_j=\{\rho:\;\arg\rho\in(\theta_j,\theta_{j+1})\},$ where
$\theta_1=\theta_5=-\varphi_2, \theta_2=\pi-\varphi_1,
\theta_3=\pi-\varphi_2, \theta_4=-\varphi_1$. For sufficiently
small $\delta>0$ we construct the sectors $S_{j,\delta}:=\{\rho:\;
\arg\rho\in(\theta_j+\delta,\theta_{j+1}-\delta)\}.$

Let $\{e_k(x,\rho)\}_{k=1,2},\; x\in[0,b]$ and
$\{E_k(x,\rho)\}_{k=1,2},\; x\in[b,T]$ be the Birkhoff-type fundamental
systems of solutions (FSS's) of Eq. (1.1) with the asymptotics as
$|\rho|\to\infty,$ $\rho\in\overline{S_j},$ $\nu=0,1$ (see [1]):
$$
e^{(\nu)}_k(x,\rho)=((-1)^{k-1}i\rho a_1)^\nu
\exp((-1)^{k-1}i\rho a_1 x)[1], \quad x\in[0,b],
$$
$$
E^{(\nu)}_k(x,\rho)=((-1)^{k-1}i\rho a_2)^\nu
\exp((-1)^{k-1}i\rho a_2(x-b))[1],\quad x\in[b,T].
$$
where $[1]=1+O(1/\rho).$
The functions $e^{(\nu)}_k(x,\rho)$ and $E^{(\nu)}_k(x,\rho)$ are
regular for $\rho\in S_j$, $|\rho|>\rho^*$ and continuous for
$\rho\in \overline{S_j}$, $|\rho|\ge\rho^*$ for some $\rho^*>0.$
Using these FSS's and the jump conditions (1.3) we get the following
asymptotical formulas as $|\rho|\to\infty,$ $\nu=0,1$:
$$
\varphi^{(\nu)}(x,\lambda)=\Big((i\rho a_1)^\nu\exp(i\rho a_1x)[1]+
(-i\rho a_1)^\nu\exp(-i\rho a_1x)[1]\Big)/2,\; x\in[0,b],
$$
$$
\varphi^{(\nu)}(x,\lambda)=\Big(\Big(\omega_{+}\exp(i\rho a_1l_1)[1]
+\omega_{-}\exp(-i\rho a_1l_1)[1]\Big)
(i\rho a_2)^\nu\exp(i\rho a_2(x-b))[1]+
$$
$$
\Big(\omega_{-}\exp(i\rho a_1l_1)[1]+
\omega_{+}\exp(-i\rho a_1l_1)[1]\Big)
(-i\rho a_2)^\nu\exp(-i\rho a_2(x-b))[1]\Big)/(4a_2),\; x\in[b,T],
$$
$$
\psi^{(\nu)}(x,\lambda)=\Big(\Big(\omega_{+}\exp(i\rho a_2l_2)[1]
-\omega_{-}\exp(-i\rho a_2l_2)[1]\Big)
(i\rho a_1)^\nu\exp(i\rho a_1(b_1-x))[1]+
$$
$$
\Big(-\omega_{-}\exp(i\rho a_2l_2)[1]
+\omega_{+}\exp(-i\rho a_2l_2)[1]\Big)
(-i\rho a_1)^\nu\exp(-i\rho a_1(b_1-x))[1]\Big)/(4a_1),\; x\in[0,b],
$$
$$
\psi^{(\nu)}(x,\lambda)=\Big((-i\rho a_2)^\nu\exp(i\rho a_2(T-x))[1]+
(i\rho a_2)^\nu\exp(-i\rho a_2(T-x))[1]\Big)/2,\; x\in[b,T].
$$
In view of (2.1), these formulas yield
$$
\Delta(\lambda)=(-i\rho)\Big(\Big(\omega_{+}\exp(i\rho a_1l_1)[1]
+\omega_{-}\exp(-i\rho a_1l_1)[1]\Big)\exp(i\rho a_2l_2)[1]-
$$
$$
\Big(\omega_{-}\exp(i\rho a_1l_1)[1]
+\omega_{+}\exp(-i\rho a_1l_1)[1]\Big)
\exp(-i\rho a_2l_2)[1]\Big)/4,\; |\rho|\to\infty,               \eqno(2.4)
$$
$$
M(\lambda)=\pm (i\rho a_1)^{-1}[1],\quad \rho\in\Pi_1^{\pm}.   \eqno(2.5)
$$
Using (2.4) by the known technique (see [1, Ch.1]) we obtain that
the spectrum $\Lambda$ of $L$ consists of two subsequences
$\Lambda=\{\lambda_{k}\}=\{\lambda_{k1}\}\cup\{\lambda_{k2}\},$ and
$$
\rho_{kj}=\sqrt{\lambda_{kj}}=\frac{k\pi}{r_jl_j}
\exp(i\theta_{3-j})+C_j+O(1/k),\; k\to\infty,                  \eqno(2.6)
$$
where $C_1=-(2ia_1l_1)^{-1}\ln(-\omega_{-}/\omega_{+}),$
$C_2=(2ia_2l_2)^{-1}\ln(\omega_{+}/\omega_{-}).$ Moreover,
$$
|\Delta(\lambda)|\ge C|\rho{\cal E}_1(\rho l_1){\cal E}_2(\rho l_2)|,
\; |M(\lambda)|\le C/|\rho|,\; \lambda\in G_{\delta},         \eqno(2.7)
$$
$$
|\varphi(x,\lambda)|\le C|{\cal E}_1(\rho x)|,
\; x\in(0,b),\quad \forall\,\lambda,
$$
$$
|\varphi(x,\lambda)|\le C|{\cal E}_1(\rho x){\cal E}_2(\rho(x-b))|,
\; x\in(b,T),\quad \forall\,\lambda,
$$
$$
|\Phi(x,\lambda)|\le C|\rho {\cal E}_1(\rho x)|^{-1},
\; x\in(0,b),\quad \lambda\in G_{\delta},
$$
$$
|\Phi(x,\lambda)|\le C|\rho {\cal E}_1(\rho x){\cal E}_2(\rho(x-b))|^{-1},
\; x\in(b,T),\quad \lambda\in G_{\delta},
$$
where $G_{\delta}:=\{\rho:\;|\rho-\rho_k|\}\ge\delta,$
${\cal E}_k(\rho x):=\exp(\pm i\rho a_kx)$ for $\rho\in\Pi_k^{\pm},$
$x\in l_k$. Let $m_k$ be the multiplicity of the eigenvalue
$\lambda_k$ ($\lambda_k=\lambda_{k+1}=\ldots=\lambda_{k+m_k-1}$), and
put $S:=\{k\ge 1:\; \lambda_{k-1}\ne\lambda_k\}\cup\{0\}.$ It follows
from (2.6) that for sufficiently large $k$ ($k>k^*$) all eigenvalues
are simple, i.e. $m_k=1$ for $k>k^*.$ Similar to [18] one gets
$$
M(\lambda)=\sum_{k\in S}\sum_{\nu=0}^{m_k-1}
\frac{M_{k+\nu}}{(\lambda-\lambda_k)^{\nu+1}}\,,                 \eqno(2.8)
$$
where $\displaystyle\sum_{\nu}
\displaystyle\frac{M_{k+\nu}}{(\lambda-\lambda_k)^{\nu+1}}$
is the principal part of $M(\lambda)$ in a neighborhood of $\lambda_k$.
The sequence ${\cal M}=\{M_k\}_{k\ge 0}$ is called the Weyl
sequence of $L,$ and the data $W=\{\lambda_k,M_k\}_{k\ge 0}$ are
called the spectral data of $L.$ Similar to (2.6) we calculate
${\cal M}=\{M_{k1}\}\cup\{M_{k2}\},$ and
$$
M_{k1}=\frac{2}{a_1^2l_1}\Big(1+O\Big(\frac{1}{k}\Big)\Big),
\; M_{k2}=\frac{8}{\omega_{-}\omega_{+}l_2}
\exp\Big(\frac{2k\pi r_1l_1}{r_2l_2}(\cos\alpha+i\sin\alpha)\Big)
\Big(1+O\Big(\frac{1}{k}\Big)\Big),                            \eqno(2.9)
$$
as $k\to\infty.$ Here $\alpha:=\varphi_1-\varphi_2+\pi/2.$
Note that $\cos\alpha<0,$ since $\alpha\in(\pi/2,3\pi/2).$
Using (2.4), (2.6), (2.7), (2.9) and the asymptotical formulas
for $\varphi(x,\lambda)$ and $\psi(x,\lambda),$ we obtain the estimates
$$
|\varphi(x,\lambda_{k1})|\le C,\quad |\varphi(x,\lambda_{k2})|
\le C\exp\Big(\frac{-k\pi r_1l_1\cos\alpha}{r_2l_2}\Big),
\quad x\in[0,T].                                              \eqno(2.10)
$$
It follows from (2.5) and (2.6) that
$$
a_1=\lim_{|\rho|\to\infty}
(i\rho M(\lambda))^{-1},\quad \rho\in\Pi_1^{+},               \eqno(2.11)
$$
$$
l_1=b=-\lim_{k\to\infty} (k\pi/(a_1\rho_{k1})),\quad l_2=T-l_1,
\quad a_2=\lim_{k\to\infty} (k\pi/(l_2\rho_{k2}),             \eqno(2.12)
$$
$$
A:=\omega_{+}/\omega_{-}=\lim_{k\to\infty}\exp(2i\rho_{k2}a_2l_2),
\quad d_1=\sqrt{(a_1(A+1))/(a_2(A-1))}.                       \eqno(2.13)
$$

\medskip
{\bf 3. Inverse problem}

\medskip
In this paper we consider the following inverse problem.

\smallskip
{\bf Inverse problem 3.1. } {\it Given the Weyl function
$M(\lambda)$ (or the spectral data $W$), construct $L.$}

\smallskip
According to (2.8) the specification of the Weyl function is
equivalent to the specification of the spectral data.

Firstly, we will prove the uniqueness theorem. For this purpose,
together with $L$ we consider a boundary value problem $\tilde L$
of the same form but with $\tilde q(x),\tilde b,\tilde r(x),
\tilde h,\tilde H,\tilde d_1,\tilde d_2$ instead of
$q(x),b,r(x),h,H,d_1,d_2$. We agree that if a certain symbol
$\chi$ denotes an object related to $L,$ then $\tilde\chi$
will denote an analogous object related to $\tilde L.$

\smallskip
{\bf Theorem 3.1. }{\it If $M(\lambda)\equiv \tilde M(\lambda)$
(or $W=\tilde W$), then $L=\tilde L.$ Thus, the specification of the
Weyl function (or the spectral data) uniquely determines the
functions $q(x), r(x)$ and the parameters $b,h,H,d_1,d_2$.}

\smallskip
{\it Proof. } It follows form (2.11)-(2.13) that $b=\tilde b,$
$a_k=\tilde a_k,$ $d_1=\tilde d_1.$ We construct the functions
$$
{\cal P}_0=\Phi\tilde\varphi-\varphi\tilde\Phi,
\quad {\cal P}_1=\varphi\tilde\Phi'-\Phi\tilde\varphi'.        \eqno(3.1)
$$
In view of (2.3), this yields
$$
\varphi={\cal P}_1\tilde\varphi+{\cal P}_0\tilde\varphi',
\quad \Phi={\cal P}_1\tilde\Phi+{\cal P}_0\tilde\Phi',
\quad{\cal P}_1-1=\varphi(\tilde\Phi'-\Phi')-
\Phi(\tilde\varphi'-\varphi').                                 \eqno(3.2)
$$
Using (2.2), (3.1), (3.2) and the asymptotical formulas for
$\varphi$ and $\psi,$ we infer
$$
|{\cal P}_1(x,\lambda)-1|\le C/|\rho|,\quad
|{\cal P}_0(x,\lambda)|\le C/|\rho|,
\quad \rho\in G_\delta\cap\tilde G_\delta.                     \eqno(3.3)
$$
Taking (2.2), (3.1) and the assumption of the theorem into account,
we conclude that the functions ${\cal P}_k(x,\lambda)$ are entire
in $\lambda$ for each $x.$ Together with (3.3) this yields
${\cal P}_1(x,\lambda)\equiv 1,$ ${\cal P}_0(x,\lambda)\equiv 0.$
Using (3.2) we calculate $\varphi(x,\lambda)
\equiv\tilde\varphi(x,\lambda),$ $\Phi(x,\lambda)\equiv
\tilde\Phi(x,\lambda),$ hence $L=\tilde L.$ Theorem 3.1 is proved.

\smallskip
Let us go on to deriving a constructive solution of the inverse problem.
For this purpose we will use ideas of the method of spectral mappings
[2]. We will reduce our nonlinear inverse problem to the solution of
the so-called {\it main equation}, which is a linear equation in a
corresponding Banach space of sequences. We give a derivation of the
main equation, and prove its unique solvability. Using the solution of
the main equation we provide an algorithm for the solution of the
inverse problem considered.
For simplicity, in the sequel we confine ourselves to the
case when the function $\Delta(\lambda)$ has only simple zeros
(the general case requires minor technical modifications).

\smallskip
Let the Weyl function $M(\lambda)$ and the spectral data $W$ be given.
Using (2.15)-(2.17) we compute $b, a_k$ and $d_1$. Then we choose
a model boundary value problem $\tilde L$ such that $\tilde b=b,
\tilde a_k=a_k, \tilde d_1=d_1$, and arbitrary in the rest (for
example, we can take $\tilde q=0$). Let $\theta_k:=1$ if $\lambda_k=
\lambda_{k1}$, and $\theta_k:=\exp(-k\pi r_1l_1(r_2l_2)^{-1}\cos\alpha)$
if $\lambda_k=\lambda_{k2}$. Denote
$$
\xi_k:=|\rho_k-\tilde\rho_k|+|M_k-\tilde M_k|\theta_k^2,\;
z_{k0}:=\lambda_k,\; z_{k1}:=\tilde\lambda_k,\; \beta_{k0}:=M_k,\;
\beta_{k1}:=\tilde M_k.
$$
By virtue of (2.6) and (2.9) one has $\xi_k=O(1/k).$
Consider the functions
$$
\varphi_{kj}(x):=\varphi(x,z_{kj}),\;
\tilde\varphi_{kj}(x):=\tilde\varphi(x,z_{kj}),\;j=0,1,
$$
$$
B_{ni,kj}(x):=D(x,z_{ni},z_{kj})\beta_{kj},\;
\tilde B_{ni,kj}(x):=\tilde D(x,z_{ni},z_{kj})\beta_{kj},\; i,j=0,1,
$$
$$
f_{k0}(x):=(\varphi_{k0}(x)-\varphi_{k1}(x))/(\xi_k\theta_k),\;
f_{k1}(x):=\varphi_{k1}(x)/\theta_k,
$$
$$
A_{n0,k0}(x):=(B_{n0,k0}(x)-B_{n1,k0}(x))\xi_k\theta_k/(\xi_n\theta_n),
$$
$$
A_{n1,k1}(x):=(B_{n1,k0}(x)-B_{n1,k1}(x))\theta_k/\theta_n,\;
A_{n1,k0}(x):=B_{n1,k0}(x)\xi_k\theta_k/\theta_n,
$$
$$
A_{n0,k1}(x):=
(B_{n0,k0}(x)-B_{n1,k0}(x)-B_{n0,k1}(x)+B_{n1,k1}(x))
\theta_k/(\xi_n\theta_n).
$$
Similarly $\tilde f_{kj}(x)$ and $\tilde A_{ni,kj}(x)$ are defined.
Using (2.6), (2.9), (2.10) and the asymptotical formulas
for $\varphi(x,\lambda)$ we get
$$
|f_{kj}(x)|,\; |\tilde f_{kj}(x)|\le C,\quad
|A_{ni,kj}(x)|, |\tilde A_{ni,kj}(x)|\le C\xi_k(|n-k|+1)^{-1}.  \eqno(3.4)
$$
Denote by $V$ the set of indices $u=(n,i),$ where $n\ge 0,\,i=0.1.$

\smallskip
{\bf Theorem 3.2. }{\it The following relation holds
$$
\tilde f_{ni}(x)=f_{ni}(x)+
\sum_{(k,j)\in V} \tilde A_{ni,kj}(x)f_{kj}(x),\; (n,i)\in V,  \eqno(3.5)
$$
where the series converge absolutely and uniformly on $x\in[0,T]$
and $\lambda$ on compact sets.}

\smallskip
{\it Proof. } Consider the contours
$\Gamma_N:=\{\lambda:\, |\lambda|=R_N\},$ where $R_N\to\infty$
such that $\Gamma_N\subset G_{\delta}$.
%and $\Lambda\cap int\,\Gamma_N=\{\lambda_k\}_{k\le N}.$
Denote ${\cal S}_k:=\{\rho:\,Im(\rho a_k)=0\},$
${\cal S}_0:={\cal S}_1\cup{\cal S}_2$,
${\cal S}:=\{\rho:\,dist({\cal S}_0,\rho)=\delta\},$ where $\delta>0$
is such that $\Lambda\cup\tilde\Lambda\subset\,int\,{\cal S}.$
Let $\gamma$ be the image of ${\cal S}$ in the $\lambda$- plane, and
$\Gamma'_N:=\Gamma_N\cap int\,\gamma,$
$\Gamma''_N:=\Gamma_N\setminus\Gamma'_N$,
$\gamma^*_N:=\gamma\cap int\,\Gamma_N$. Denote by
$\gamma_N:=\gamma^*_N\cup\Gamma'_N$ and
$\gamma^0_N:=\gamma^*_N\cup\Gamma''_N$
the closed contours with counterclockwise circuit.
Applying Cauchy's integral formula we get
$$
{\cal P}_k(x,\lambda)-\delta_{1k}=\frac{1}{2\pi i}\int_{\gamma^0_N}
\frac{{\cal P}_k(x,\mu)-\delta_{1k}}{\lambda-\mu}d\mu=
\frac{1}{2\pi i}\int_{\gamma_N}\frac{{\cal P}_k(x,\mu)}{\lambda-\mu}d\mu
-\frac{1}{2\pi i}\int_{\Gamma_N}
\frac{{\cal P}_k(x,\mu)-\delta_{1k}}{\lambda-\mu}d\mu,
$$
where $k=0,1,$ $\lambda\in int\,\gamma^0_N$, and $\delta_{jk}$ is the
Kronecker delta. Taking (3.2) into account we calculate
$$
\varphi(x,\lambda)=\tilde\varphi(x,\lambda)+\frac{1}{2\pi i}
\int_{\gamma_N}\Big(\tilde\varphi(x,\lambda){\cal P}_1(x,\mu)+
\tilde\varphi'(x,\lambda){\cal P}_0(x,\mu)\Big)\frac{d\mu}{\lambda-\mu}
+\varepsilon_N(x,\lambda).
$$
In view of (3.3), one has $\displaystyle\lim_{N\to\infty}
\varepsilon_N(x,\lambda)=0$
uniformly in $x\in[0,T]$ and $\lambda$ on compact sets.
Taking (3.1) and (2.2) into account we obtain
$$
\tilde\varphi(x,\lambda)=\varphi(x,\lambda)+\frac{1}{2\pi i}
\int_{\gamma_N}\tilde D(x,\lambda,\mu)(M(\mu)-
\tilde M(\mu))\varphi(x,\mu)d\mu+\varepsilon_N(x,\lambda).
$$
Note that the terms with $S(x,\lambda)$ and $\tilde S(x,\lambda)$ are
zero because of Cauchy's theorem. Using the residue theorem we get the
relation
$$
\tilde\varphi_{ni}(x)=\varphi_{ni}(x)+\sum_{k=0}^\infty
\Big(\tilde B_{ni,k0}(x)\varphi_{k0}(x)
-\tilde B_{ni,k1}(x)\varphi_{k1}(x)\Big),
$$
which is equivalent to (3.5). Theorem 3.2 is proved.

\smallskip
By similar arguments we calculate
$$
A_{ni,kj}(x)-\tilde A_{ni,kj}(x)+\sum_{(l,s)\in V}
\tilde A_{ni,ls}(x)A_{ls,kj}(x)=0,\; (n,i),(k,j)\in V.        \eqno(3.6)
$$
Let $f(x)=[f_u(x)]_{u\in V}$, $A(x)=[A_{u,v}(x)]_{u,v\in V}$,
$\tilde f(x)=[\tilde f_u(x)]_{u\in V}$,
$\tilde A(x)=[\tilde A_{u,v}(x)]_{u,v\in V}$. We denote by $m$
the Banach space of bounded sequences $\chi=[\chi_u]_{u\in V}$
with the norm $\|\chi\|=\sup_{u\in V} |\chi_u|.$
According to (3.4), one has that for each fixed $x,$ the operators
$I+\tilde A(x)$ and $I-A(x),$ acting from $m$ to $m,$ are linear
bounded operators. Relations (3.5) and (3.6) can be written as follows
$$
\tilde f(x)=(I+\tilde A(x))f(x),\quad (I+\tilde A(x))(I-A(x))=I.
$$
Symmetrically one has
$f(x)=(I-A(x))\tilde f(x),\quad (I-A(x))(I+\tilde A(x))=I.$
Thus, for each fixed $x,$ the operator $I+\tilde A(x)$ has
a bounded inverse operator, hence the linear equation
$\tilde f(x)=(I+\tilde A(x))f(x)$ is uniquely solvable. This
equation is called {\it the main equation} of the inverse problem.
Solving the main equation we find the vector $f(x),$ and also the
solutions $\varphi_{ni}(x)=\varphi(x,\lambda_{ni})$ of Eq. (1.1), hence
we can construct $q(x), h, H$ and $d_2$. Thus, the solution of the
inverse problem can be found by the following algorithm.

\smallskip
{\bf Algorithm 3.1. }{\it Given the Weyl function $M(\lambda)$ and
the spectral data $W.$ \\
1) Calculate $b, a_k$ and $d_1$ via (2.11)-(2.13).\\
2) Choose a model boundary value problem $\tilde L$ such that
$\tilde b=b, \tilde a_k=a_k, \tilde d_1=d_1$.\\
3) Construct $\tilde f(x)$ and $\tilde A(x)$ (see above).\\
4) Find $f(x)=[f_u]_{u\in V}$ by solving the main equation
$\tilde f(x)=(I+\tilde A(x))f(x).$\\
5) Calculate $\varphi_{n1}(x)=f_{n1}(x)\theta_n,
\varphi_{n0}=\varphi_{n1}(x)+f_{n0}(x)\xi_n\theta_n$.\\
6) Find $q(x), h, H$ and $d_2$ using (1.1)-(1.3).}

\smallskip
{\bf Remark 3.1. } We can also calculate $q(x)$ by the formula
$q(x)=\tilde q(x)-2{\cal F}(x),$ where
$$
{\cal F}(x)=\frac{d}{dx}
\sum_{k=0}^{\infty}\Big(M_{k0}\tilde\varphi_{k0}(x)\varphi_{k0}(x)
-M_{k1}\tilde\varphi_{k1}(x)\varphi_{k1}(x)\Big).              \eqno(3.7)
$$

\begin{center}
{\bf REFERENCES}
\end{center}
\begin{enumerate}
\item[{[1]}] {\it Freiling G. and Yurko V.A.} Inverse Sturm-Liouville
       Problems and their Applications. NOVA Science Publishers.  New York, 2001.
\item[{[2]}] {\it Yurko V.A.} Method of Spectral Mappings in the Inverse
       Problem Theory. Inverse and Ill-posed Problems Series. VSP. Utrecht,  2002.
\item[{[3]}] {\it Yurko V.A.}  Introduction to the Theory of Inverse
     Spectral Problems. Fizmatlit, Moscow, 2007 [in Russian].
\item[{[4]}] {\it Yurko V.A.} Inverse spectral problems for differential
      operators on spatial  networks. Russian Mathematical Surveys,
      vol.71, no.3 (2016), 539-584.
\item[{[5]}]  {\it Krueger R.J.} Inverse problems for nonabsorbing media
     with discontinuous material properties. J. Math.Phys. 23, no.3 (1982),
     396–404.
\item[{[6]}]  {\it Anderssen R.S.} The effect of discontinuities in density
      and shear velocity on the asymptotic overtone structure ofrtional
      eigenfrequencies of the Earth. Geophys. J. Int. 50, no.2 (1977), 303–309.
\item[{[7]}] {\it Hald O.H.} Discontinuous inverse eigenvalue problems.
      Comm. Pure Appl.Math. 37, no.5 (1984), 539–577.
\item[{[8]}] {\it Yurko V.A.} Boundary-value problems with discontinuity
      conditions at an interior point of the interval. Differ. Uravn.
      36, no.8 (2000), 1139–1140. English transl. in Differ. Equations
      36, no. 8 (2000), 1266–1269.
\item[{[9]}] {\it Yurko V.A.} Integral transforms connected with
     discontinuous boundary value problems. Integral Transform. Spec.
     Functions, 10, no.2 (2000), 141–164.
\item[{[10]}] {\it Belishev M.} An inverse spectral indefinite problem
      for the equation $y''+zp(x)y=0$ on an interval. Funct. Anal. Appl.
      21, no.2 (1987), 68–69.
\item[{[11]}] {\it Daho K. and Langer H.}  Sturm–Liouville operators with
     an indefinite weight  functions. Proc. R. Soc. Edinb. Sect. A 78
     (1977), 161–191.
\item[{[12]}]  {\it Andersson L.-E.} Inverse eigenvalue problems with
      discontinuous coefficients. Inverse Problems 4, no.2 (1988), 353–397.
\item[{[13]}] {\it Coleman C. and McLaughlin J.} Solution of the inverse
      spectral problem for an impedance with integrable derivative, I, II.
      Comm. Pure Appl. Math. 46, no.2, (1993) 145–184; Comm. Pure Appl.
      Math. 46, no.2 (1993), 185-212.
\item[{[14]}] {\it Freilng G. and Yurko V.} Inverse problems for
     differential equations with turning points.  Inverse Problems, 13,
     no.5 (1997), 1247–1263.
\item[{[15]}] {\it Yurko V.A.} Inverse spectral problems for Sturm-Liouville
      operators with complex weights. Inverse Problems in Science and
      Engineering, 26, no.10 (2018), 1396-1403.
\item[{[16]}] {\it Yurko V.A.} An inverse problem for Sturm-Liouville
      operators on the half-line with complex weights. Journal of Inverse
      and Ill-Posed Problems, 27, no.3 (2019),  439–443.
\item[{[17]}] {\it Golubkov A.A. and Kuryshova Yu.V.} Inverse problem for
     Sturm–Liouville operators on a curve. Tamkang Journal of Mathematics,
     50, no.3 (2019), 349-359.
\item[{[18]}] {\it Buterin S.} On inverse spectral problems for
     non-selfadjoint St-L operator on a finite interval. J. Math.
     Analysis and Appl. 335 (2007), 739-749.\\
     https://doi.org/10.1016/j.jmaa.2007.02.012
\end{enumerate}

\noindent {\bf Vjacheslav A. Yurko}, yurkova@info.sgu.ru,
https://orcid.org/0000-0002-4853-0102

\end{document}